\documentclass[preprint,showpacs,preprintnumbers,amsmath,amssymb]{revtex4}
\input{epsf.sty}
\usepackage{epsfig}
\usepackage{graphicx}
\usepackage{dcolumn}
\usepackage{bm}
\usepackage{times}

\usepackage{amsmath}

\usepackage{amssymb}

\newtheorem{theorem}{Theorem}
\newtheorem{lemma}{Lemma}

\newcommand{\onetom}{1,\cdots,m}

\newcommand{\onetoK}{1,\cdots,K}

\begin{document}
\title{Cluster synchronization in networks of coupled non-identical dynamical systems}

\author{Wenlian Lu}
\email{wenlian@fudan.edu.cn} \affiliation{Centre for Computational
Systems Biology, Fudan University, Shanghai, P. R. China}
\author{Bo Liu}
\email{071018024@fudan.edu.cn}
\author{Tianping Chen}
\email{tchen@fudan.edu.cn} \affiliation{Shanghai Key Laboratory for
Contemporary Applied Mathematics, School of Mathematical Sciences,
Fudan University, Shanghai, P. R. China}

\date{\today}

\begin{abstract}
In this paper, we study cluster synchronization in networks of
coupled non-identical dynamical systems. The vertices in the same
cluster have the same dynamics of uncoupled node system but the
uncoupled node systems in different clusters are different. We
present conditions guaranteeing cluster synchronization and
investigate the relation between cluster synchronization and the
unweighted graph topology. We indicate that two condition play key
roles for cluster synchronization: the common inter-cluster coupling
condition and the intra-cluster communication. From the latter one,
we interpret the two well-known cluster synchronization schemes:
self-organization and driving, by whether the edges of communication
paths lie at inter- or intra-cluster. By this way, we classify
clusters according to whether the set of edges inter- or
intra-cluster edges are removable if wanting to keep the
communication between pairs of vertices in the same cluster. Also,
we propose adaptive feedback algorithms on the weights of the
underlying graph, which can synchronize any bi-directed networks
satisfying the two conditions above. We also give several numerical
examples to illustrate the theoretical results.

\end{abstract}


\maketitle {\bf Cluster synchronization is considered to be more
momentous than complete synchronization in brain science and
engineering control, ecological science and communication
engineering, social science and distributed computation. Most of the
existing works only focused on networks with either special
topologies such as regular lattices or coupled two/three groups. For
the general coupled dynamical systems, theoretical analysis to
clarify the relationship between the (unweighted) graph topology and
the cluster scheme, including both self-organization and driving, is
absent. In this paper, we study this topic and find two essential
conditions for an unweighted graph topology to realize cluster
synchronization: the common inter-cluster coupling condition and the
intra-cluster communication. Thus under these conditions, we present
two manners of weighting to achieve cluster synchronization. One is
adding positive weights on each edges with keeping the invariance of
the cluster synchronization manifold and the other is an adaptive
feedback weighting algorithms. We prove the availability of each
manner. From these results, we give an interpretation of the two
clustering synchronization schemes: self-organization and driving,
involved with the unweighted graph topology, via the communication
between pairs of individuals in the same cluster. Thus, we present
one way to classify the clusters via whether the set of inter- or
intra-cluster edges are removable if still wanting to keep the
communication between vertices in the same cluster. }
\section{Introduction}
Recent decades witnesses that chaos synchronization in complex
networks has attracted increasing interests from many research and
application fields \cite{Pik,Boc1,Wang1}, since it was firstly
introduced in Ref. \cite{Yamada}. Word ``synchronization'' comes
from Greek, which means ``share time'' and today, it comes to be
considered as ``time coherence of different processes''. Many new
synchronization phenomena appear in a wide range of real systems,
such as biology \cite{Str}, neural networks \cite{Gra},
physiological processes \cite{Gla}. Among them, the most interesting
cases are complete synchronization, cluster synchronization, phase
synchronization, imperfect synchronization, lag synchronization, and
almost synchronization etc. See Ref. \cite{Boc} and the references
therein.

Complete synchronization is the most special one and characterized
by that all oscillators approach to a uniform dynamical behavior. In
this situation, powerful mathematical techniques from dynamical
systems and graph theory can be utilized. Pecora et.al. \cite{Pec}
proposed the Master Stability Function for transverse stability
analysis \cite{Ash} of the diagonal synchronization manifold. This
method has been widely used to study local completer synchronization
in networks of coupled system \cite{Loc_syn}. Refs.
\cite{Wu,Glo_syn,Lu} proposed a framework of Lyapunov function
method to investigate global synchronization in complex networks.
One of the most important issues is how the graph topology affects
the synchronous motion \cite{Boc1}. As pointed out in Ref.
\cite{Wu2}, the connectivity of the graph plays a significant role
for chaos synchronization.

Cluster synchronization is considered to be more momentous in brain
science \cite{Sch} and engineering control \cite{Clu_con},
ecological science \cite{Clu_eco} and communication engineering
\cite{Clu_com}, social science \cite{Sto} and distributed
computation \cite{Clu_dis}. This phenomenon is observed when the
oscillators in networks are divided into several groups, called
clusters, by the way that all individuals in the same cluster reach
complete synchronization but the motions in different clusters do
not coincide. Cluster synchronization of coupled identical systems
are studied in Refs. \cite{Bel1,Liu,Wei,Jalan}. Among them, Jalan
et. al. \cite{Jalan} pointed out two basic formations which realize
cluster synchronization. One is {\em self-organization}, which leads
to cluster with dominant intra-cluster couplings, and the other is
{\em driving}, which leads to cluster with dominant inter-cluster
couplings.

Nowadays, the interest of cluster synchronization is shifting to
networks of coupled non-identical dynamical systems. In this case,
cluster synchronization is obtained via two aspects: the oscillators
in the same cluster have the same uncoupled node dynamics and the
inter- or intra-cluster interactions realize cluster synchronization
via driving or/and self-organizing configurations. Refs. \cite{Liu}
proposed cluster synchronization scheme via dominant intra-couplings
and common inter-cluster couplings. Ref. \cite{Sor} studied local
cluster synchronization for bipartite systems, where no
intra-cluster couplings (driving scheme) exist. Refs. \cite{ChenL}
investigated global cluster synchronization in networks of two
clusters with inter- and intra-cluster couplings. Belykh et. al.
studied this problem in 1D and 2D lattices of coupled identical
dynamical systems in Ref. \cite{Bel1} and non-identical dynamical
systems in Ref. \cite{Bel2}, where the oscillators are coupled via
inter- or/and intra-cluster manners. Ref. \cite{Pha} used nonlinear
contraction theory \cite{Loh} to build up a sufficient condition for
the stability of certain invariant subspace, which can be utilized
to analyze cluster synchronization (concurrent synchronization is
called in that literature). However, up till now, there are no works
revealing the relationship between the (unweighted) graph topology
and the cluster scheme, including both self-organization and
driving, for general coupled dynamical systems.

The purpose of this paper is to study cluster synchronization in
networks of coupled non-identical dynamical systems with various
graph topologies. In Section 2, we formulate this problem and study
the existence of the cluster synchronization manifold. Then, we give
one way to set positive weights on each edges and derive a criterion
for cluster synchronization. This criterion implies that the
communicability between each pair of individuals in the same cluster
is essential for cluster synchronization. Thus, we interpret the two
communication schemes: self-organization and driving, according to
the communication scheme among individuals in the same cluster. By
this way, we classify clusters according to the manner by which
synchronization in a cluster realizes. In Sec. 3, we propose an
adaptive feedback algorithms on weights of the graph to achieve a
given clustering. In Sec. 5, we discuss the cluster
synchronizability of a graph with respect to a given clustering and
present the general results for cluster synchronization in networks
with general positive weights. We conclude this paper in Sec. 6.

\section{Cluster synchronization analysis}
In this section, we study cluster synchronization in a network with
weighted bi-directed graph and a given division of clusters. We
impose the constraints on graph topology to guarantee the invariance
of the corresponding cluster synchronization manifold and derive the
conditions for this invariant manifold to be globally asymptotically
stable by the Lyapunov function method. Before that, we should
formulate the problem.

Throughout the paper, we denote a positive definite matrix $Z$ by
$Z>0$ and similarly for $Z<0$, $Z\le 0$, and $Z\ge 0$. We say that a
matrix $Z$ is positive definite on a linear subspace $V$ if
$u^{\top} Z u>0$ for all $u\in V$ and $u\ne 0$, denoted by
$Z|_{V}>0$. Similarly, we can define $Z|_{V}<0$, $Z|_{V}\ge 0$, and
$Z|_{V}\le 0$. If a matrix $Z$ has all eigenvalues real, then we
denote by $\lambda_{k}(Z)$ the $k$-th largest eigenvalues of $Z$.
$Z^{\top}$ denotes the transpose of the matrix $Z$ and
$Z^{s}=(Z+Z^{\top})/2$ denotes the symmetry part of a square matrix
$Z$. $\# A$ denotes the number of the set $A$ with finite elements.
\subsection{Model description and existence of invariant cluster synchronization
manifold}

A bi-directed unweighted graph $\mathcal G$ is denoted by a double
set $\{\mathcal V,\mathcal E\}$, where $\mathcal V$ is the vertex
set numbered by $\{\onetom\}$, and $\mathcal E$ denotes the edge set
with $e(i,j)\in\mathcal E$ if and only if there is an edge
connecting vertices $j$ and $i$. $\mathcal N(i)=\{j\in\mathcal
V:~e(i,j)\in\mathcal E\}$ denotes the neighborhood set of vertex
$i$. The graph considered in this paper is always supposed to be
simple (without self-loops and multiple edges) and bi-directed. A
clustering $\mathcal C$ is a disjoint division of the vertex set
$\mathcal V$: $\mathcal C=\{\mathcal C_{1},\mathcal
C_{2},\cdots,\mathcal C_{K}\}$ satisfying (i).
$\bigcup_{k=1}^{K}\mathcal C_{k}=\mathcal V$; (ii). $\mathcal
C_{k}\bigcap\mathcal C_{l}=\emptyset$ holds for $k\ne l$.

The network of coupled dynamical system is defined on the graph
$\mathcal G$. The individual uncoupled system on the vertex $i$ is
denoted by an $n$-dimensional ordinary differential equation
$\dot{x}^{i}=f_{k}(x^{i})$ for all $i\in\mathcal C_{k}$, where
$x^{i}=[x^{i}_{1},\cdots,x^{i}_{n}]^{\top}$ is the state variable
vector on vertex $i$ and $f_{k}(\cdot):\mathbb R^{n}\to\mathbb
R^{n}$ is a continuous vector-valued function. Each vertex in the
same cluster has the same individual node dynamics. The interaction
among vertices is denoted by linear diffusion terms. It should be
emphasized that $f_{k}$ for different clusters are distinct, which
can guarantee that the trajectories are apparently distinguishing
when cluster synchronization is reached.

Consider the following model of networks of linearly coupled
dynamical system \cite{IEEE}:
\begin{eqnarray}
\dot{x}^{i}=f_{k}(x^{i})+\sum_{j\in\mathcal
N(i)}w_{ij}\Gamma(x^{j}-x^{i}),~i\in\mathcal
C_{k},~k=\onetoK.\label{Eq.1}
\end{eqnarray}
where $w_{ij}$ is the coupling weight at the edge from vertex $j$ to
$i$ and $\Gamma=[\gamma_{uv}]_{u,v=1}^{n}$ denotes the inner
connection by the way that $\gamma_{uv}\ne 0$ if the the $u$-th
component of the vertices can be influenced by the $v$-th component.
The graph $\mathcal G$ is bi-directed and the weights are not
requested to be symmetric. Namely, we don't request $w_{ij}= w_{ji}$
for each pair $(i,j)$ with $e(i,j)\in\mathcal E$.

Let $A=[a_{ij}]_{i,j=1}^{m}$ be the adjacent matrix of the graph
$\mathcal G$. That is, $a_{ij}=1$ if $e(i,j)\in\mathcal E$;
$a_{ij}=0$ otherwise. Then, model (\ref{Eq.1}) can be rewritten as
\begin{eqnarray}
\dot{x}^{i}=f_{k}(x^{i})+\sum_{j=1}^{m}a_{ij}w_{ij}\Gamma(
x^{j}-x^{i}),~i\in\mathcal C_{k},~k=\onetoK.\label{Eq.2}
\end{eqnarray}
In this paper, cluster synchronization is defined as follows:
\begin{enumerate}
\item The
differences among trajectories of vertices in the same cluster
converge to zero as time goes to infinity, i.e.,
\begin{eqnarray}
\lim_{t\to\infty}[x^{i}(t)-x^{j}(t)]=0,~\forall~i,j\in\mathcal
C_{k},~k=\onetoK;\label{cs}
\end{eqnarray}
\item The differences among the trajectories of vertices in
different clusters do not converge to zero, i.e.,
$\overline{\lim}_{t\to\infty}|x^{i'}(t)-x^{j'}(t)|>0$ holds for each
$i'\in\mathcal C_{k}$ and $j'\in\mathcal C_{l}$ with $k\ne l$.
\end{enumerate}
As mentioned above, we suppose that the latter one can be guaranteed
by the incoincidence of $f_{k}(\cdot)$. Under this prerequisite
assumption, cluster synchronization is equivalent to the
asymptotical stability of the following cluster synchronization
manifold with respect to the clustering $\mathcal C$:
\begin{eqnarray}
\mathcal S_{\mathcal
C}(n)=\{[{x^{1}}^{\top},\cdots,{x^{m}}^{\top}]^{\top}:~x^{i}=x^{j}\in\mathbb
R^{n},~\forall~i,j\in\mathcal C_{k},~k=\onetoK\}.
\end{eqnarray}

To investigate cluster synchronization, a prerequisite requirement
is that the manifold $\mathcal S_{\mathcal C}(n)$ should be
invariant through Eqs. (\ref{Eq.2}). Assume that $x^{i}(t)=s^{k}(t)$
for each $i\in\mathcal C_{k}$ is the synchronized solution of the
cluster $\mathcal C_{k}$, $k=\onetoK$. By Eqs. (\ref{Eq.2}), each
$s^{k}$ must satisfy
\begin{eqnarray}
\dot{s}^{k}=f_{k}(s^{k})+\sum_{k'=1,k'\ne
k}^{K}\alpha_{i,k'}\Gamma(s^{k'}-s^{k}), ~~\forall~i\in\mathcal
C_{k}, \label{cluster_syn}
\end{eqnarray}
where $\alpha_{i,k'}=\sum_{j\in\mathcal C_{k'}}a_{ij}w_{ij}$.  This
demands $\alpha_{i_{1},k'}=\alpha_{i_{2},k'}$ for any $i_{1}\in
\mathcal C_{k}$, $i_{2}\in \mathcal C_{k}$, namely, $\alpha_{i,k'}$
is independent of $i$. Therefore, we have
\begin{eqnarray}
\alpha_{i,k'}=\alpha(k,k'),~i\in\mathcal C_{k},~k\ne
k'.\label{common}
\end{eqnarray}
This condition is sufficient and necessary for  the cluster
synchronization manifold $\mathcal S_{\mathcal C}(n)$ is invariant
through the coupled system (\ref{Eq.2}) for general maps
$f_{k}(\cdot)$.

Denote $\mathcal N_{k'}(i)=\mathcal N(i)\bigcap\mathcal C_{k'}$, and
define an index set $\mathcal L^{i}_{k}=\{k':~k'\ne k,~{\rm
and}~N_{k'}(i)\ne\emptyset\}$. The set $\mathcal L^{i}_{k}$
represents those clusters other than $\mathcal{C}_{k}$ and have
links to the vertex $i$. To satisfy the condition (\ref{common}),
the following {\em common inter-cluster coupling condition} over the
unweighted graph topology should be satisfied: for $k=\onetoK$,
\begin{eqnarray}
\mathcal L^{i}_{k}=\mathcal L^{i'}_{k},~\forall~i,i'\in\mathcal
C_{k}.\label{inter-cluster}
\end{eqnarray}
Therefore, we can use $\mathcal L_{k}$ to represent $\mathcal
L_{k}^{i}$ for all $i\in\mathcal C_{k}$ if the common inter-cluster
coupling condition is satisfied.

Throughout this paper, we assume that the vector-valued function
$f_{k}(x)-\alpha\Gamma x:\mathbb R^{n}\to\mathbb R^{n}$ satisfies
{\em decreasing condition} for some $\alpha\in\mathbb R$. That is,
there exists $\delta>0$ such that
\begin{eqnarray}
(\xi-\zeta)^{\top}\bigg[f_{k}(\xi)-f_{k}(\zeta)-\alpha\Gamma(\xi-\zeta)\bigg]\le-\delta
(\xi-\zeta)^{\top}(\xi-\zeta).\label{decreasing}
\end{eqnarray}
holds for all $\xi,\zeta\in\mathbb R^{n}$.  This condition holds for
any globally Lipschitz continuous function $f(\cdot)$ for
sufficiently large $\alpha>0$ and $\Gamma=I_{n}$. However, even
though $f(\cdot)$ is only locally Lipschitz, if the solution of the
coupled system (\ref{Eq.1}) is essentially bounded, then restricted
to such bounded region, the condition (\ref{decreasing}) also holds
for sufficiently large $\alpha$ and $\Gamma=I_{n}$. In this paper,
we suppose that the solution of the coupled system (\ref{Eq.2}) is
essentially bounded.

\subsection{Cluster synchronization analysis}
In the following, we investigate cluster synchronization of networks
of coupled non-identical dynamical systems withe the following
weighting scheme:
\begin{eqnarray}
w_{ij}=\left\{\begin{array}{ll}\frac{c}{d_{i,k}}&j\in\mathcal
N_{k}(i)~{\rm and}~\mathcal N_{k}(i)\in\emptyset\\
0&{\rm otherwise},\end{array}\right.
\end{eqnarray}
where $d_{i,k'}=\#\mathcal N_{k}(i)$ denotes the number of elements
in $N_{k'}(i)$ and $c$ denotes the coupling strength.. Thus, the
coupled system becomes:
\begin{eqnarray}
\dot{x}^{i}=f_{k}(x^{i})
+c\bigg[\sum_{\mathcal{N}_{k'}(i)\ne\emptyset}\frac{1}{d_{i,k'}}
\sum_{j\in\mathcal N_{k'}(i)}\Gamma(x^{j}-x^{i})\bigg],~i\in\mathcal
C_{k},~k=\onetoK.\label{Eq.3}
\end{eqnarray}
It can be seen that in Eqs. (\ref{Eq.3}), for each $i\in\mathcal
C_{k}$, the corresponding $\alpha_{i,k'}=c$ for all $k'\in\mathcal
L_{k}$ under the common inter-cluster coupling condition. The
general situation can be handled by the same approach and will be
presented in the discussion section.

We denote the weighted Laplacian of the graph as follows. For each
pair $(i,j)$ with $i\ne j$, $l_{ij}=\frac{1}{d_{i,k}}$ if
$j\in\mathcal N_{k}(i)$ and $\mathcal N_{k}(i)\ne\emptyset$ for some
$k\in\{\onetoK\}$, and $l_{ij}=0$ otherwise;
$l_{ii}=-\sum_{j=1}^{m}l_{ij}$. Thus, Eqs. (\ref{Eq.3}) can be
rewritten as:
\begin{eqnarray}
\dot{x}^{i}=f_{k}(x^{i})+c\sum_{j=1}^{m}l_{ij}\Gamma
x^{j},~i\in\mathcal C_{k},~k=\onetoK.\label{Eq.4}
\end{eqnarray}

The approach to analyze cluster synchronization is extended from
that used in Ref. \cite{Lu} to study complete synchronization. Let
$d=[d_{1},\cdots,d_{m}]^{\top}$ be a vector with $d_{i}>0$ for all
$i=\onetom$. We use the vector $d$ to construct a (skew) projection
of $x=[{x^{1}}^{\top},\cdots,{x^{m}}^{\top}]^{\top}$ onto the
cluster synchronization manifold $\mathcal S_{\mathcal C}(n)$.
Define an average state with respect to $d$ in the cluster $\mathcal
C_{k}$ as
\begin{eqnarray*}
\bar{x}^{k}_{d}=\frac{1}{\sum_{i\in\mathcal
C_{k}}d_{i}}\sum_{i\in\mathcal C_{k}}d_{i}x^{i}.
\end{eqnarray*}
Thus, we denote the projection of $x$ on the cluster synchronization
manifold $\mathcal S_{\mathcal C}(n)$ with respect to $d$ as:
$\bar{x}_{d}=[\tilde{x}^{1^{\top}},\cdots,\tilde{x}^{m^{\top}}]^{\top}$
is denoted as:
\begin{eqnarray*}
\tilde{x}^{i}=\bar{x}_{d}^{k},~if~i\in\mathcal C_{k}.
\end{eqnarray*}

Then, the variations $x^{i}-\bar{x}^{k}_{d}$ compose the transverse
space:
\begin{eqnarray*}
\mathcal T_{\mathcal
C}^{d}(n)=\bigg\{u=[{u^{1}}^{\top},\cdots,{u^{m}}^{\top}]^{\top}\in\mathbb
R^{mn}:~u^{i}\in\mathbb R^{n},~\sum_{i\in\mathcal
C_{k}}d_{i}u^{i}=0,~\forall~k=\onetoK\bigg\}.
\end{eqnarray*}
In particular, in the case of $n=1$, it denotes
\begin{eqnarray*}
\mathcal T_{\mathcal
C}^{d}(1)=\bigg\{u=[u^{1},\cdots,u^{m}]^{\top}\in\mathbb
R^{m}:~\sum_{i\in\mathcal
C_{k}}d_{i}u^{i}=0,~\forall~k=\onetoK\bigg\}.
\end{eqnarray*}

From the definition, we have the following lemma which is repeatedly
used below.
\begin{lemma} \label{lem1}For each $k\in\onetoK$, it holds
\begin{eqnarray*}
\sum_{i\in\mathcal C_{k}}d_{i}(x^{i}-\bar{x}_{d}^{k})=0.
\end{eqnarray*}
\end{lemma}
In fact, note
\begin{eqnarray*}
\sum_{i\in\mathcal
C_{k}}d_{i}(x^{i}-\bar{x}_{d}^{k})=\sum_{i\in\mathcal
C_{k}}d_{i}x^{i}-\sum_{i\in\mathcal
C_{k}}d_{i}\bigg(\frac{1}{\sum_{j\in\mathcal
C_{k}}d_{j}}\bigg)\sum_{i'\in\mathcal
C_{k}}d_{i'}x^{i'}=\sum_{i\in\mathcal
C_{k}}d_{i}x^{i}-\sum_{i'\in\mathcal C_{k}}d_{i'}x^{i'}=0.
\end{eqnarray*}
The lemma immediately follows. As a direct consequence, we have
\begin{eqnarray*}
\sum_{i\in\mathcal
C_{k}}d_{i}(x^{i}-\bar{x}_{d}^{k})^{\top}J_{k}=\bigg[\sum_{i\in\mathcal
C_{k}}d_{i}(x^{i}-\bar{x}_{d}^{k})\bigg]^{\top}J_{k}=0
\end{eqnarray*}
for any $J_{k}$, with a proper dimension, independent of the index
$i$.

Since the dimension of $\mathcal T_{\mathcal C}^{d}(n)$ is $n(m-K)$,
the dimension of $\mathcal S_{\mathcal C}$ is $nK$, and $\mathcal
S_{\mathcal C}(n)$ is disjoint with $\mathcal T_{\mathcal C}^{d}(n)$
except the origin, $\mathbb R^{mn}=\mathcal S_{\mathcal
C}(n)\bigoplus\mathcal T_{\mathcal C}^{d}(n)$, where $\bigoplus$
denotes the direct sum of linear subspaces. With these notations,
the cluster synchronization is equivalent to the {\em transverse
stability} of the cluster synchronization manifold $\mathcal
S_{\mathcal C}(n)$, i.e.,  the projection of $x$ on the transverse
space $\mathcal T_{\mathcal C}^{d}(n)$ converges to zero as time
goes to infinity.

\begin{theorem}\label{st1}
Suppose that the common inter-cluster coupling condition
(\ref{inter-cluster}) holds, $\Gamma$ is symmetry and nonnegative
definite, and each vector-valued function $f_{k}(\cdot)-\alpha\Gamma
\cdot$ satisfies the decreasing condition (\ref{decreasing}) for
some $\alpha\in\mathbb R$. If there exists a positive definite
diagonal matrix $D$ such that the restriction of $[D(cL+\alpha
I_{m})]^{s}$, restricted to the transverse space $\mathcal
T_{\mathcal C}^{d}(1)$, is non-positive definite, i.e.,
\begin{eqnarray}
\bigg[D(cL+\alpha I_{m})\bigg]^{s}\bigg |_{\mathcal T_{\mathcal
C}^{d}(1)}\le 0\label{state1}
\end{eqnarray}
holds, then the coupled system (\ref{Eq.4}) can clustering
synchronize with respect to the clustering $\mathcal C$.
\end{theorem}
{\em Proof.} We define an auxiliary function to measure the distance
from $x$ to the cluster synchronization manifold as follows
\begin{eqnarray*}
V_{k}=\frac{1}{2}\sum_{i\in\mathcal
C_{k}}d_{i}(x^{i}-\bar{x}_{d}^{k})^{\top}(x_{i}-\bar{x}_{d}^{k}),~V(x)=\sum_{k=1}^{K}V_{k}.
\end{eqnarray*}
Differentiating $V_{k}$ along Eqs. (\ref{Eq.4}) gives
\begin{eqnarray*}
\dot{V}_{k}&=&\sum_{i\in\mathcal
C_{k}}d_{i}(x^{i}-\bar{x}^{k}_{d})^{\top}\bigg[f_{k}(x^{i})
+c\sum_{j=1}^{m}l_{ij}\Gamma x^{j}-\dot{\bar{x}}^{k}_{d}\bigg].
\end{eqnarray*}
Recalling the definitions of $l_{ij}$ and the common inter-cluster
coupling condition (\ref{inter-cluster}), we have
\begin{eqnarray}
\sum_{j\in\mathcal{C}_{k'}}l_{ij}
=\sum_{j\in\mathcal{C}_{k'}}l_{i'j},~\forall~i,i'\in\mathcal
C_{k},~k\ne k',\label{A}
\end{eqnarray}
which leads
\begin{eqnarray}
\sum_{j\in\mathcal{C}_{k}}l_{ij}
=\sum_{j\in\mathcal{C}_{k}}l_{i'j},~\forall~i,i'\in\mathcal
C_{k}.\label{B}
\end{eqnarray}

By Lemma \ref{lem1}, we have
\begin{eqnarray*}
&&\sum_{i\in\mathcal
C_{k}}d_{i}(x^{i}-\bar{x}^{k}_{d})^{\top}\dot{\bar{x}}^{k}_{d}=0,
\quad \sum_{i\in\mathcal
C_{k}}d_{i}(x^{i}-\bar{x}^{k}_{d})^{\top}f_{k}(\bar{x}^{k}_{d})=0,\\
&&\sum_{i\in\mathcal C_{k}}d_{i}(x^{i}-\bar{x}^{k}_{d})^{\top}
\bigg(\sum_{j\in\mathcal{C}_{k'}}l_{ij}\Gamma{\bar{x}}^{k'}_{d}\bigg)=0,~k'=1,\cdots,K,
\end{eqnarray*}
due to the facts (\ref{A}) and (\ref{B}). Therefore, we have
\begin{eqnarray*}
\dot{V}_{k}&=&\sum_{i\in\mathcal
C_{k}}d_{i}(x^{i}-\bar{x}^{k}_{d})^{\top}\bigg[f_{k}(x^{i})-f_{k}(\bar{x}_{d}^{k})+f_{k}(\bar{x}_{d}^{k})
+c\sum_{j=1}^{m}l_{ij}\Gamma
(x^{j}-\bar{x}_{d}^{k'})-\dot{\bar{x}}^{k}_{d}
+c\sum_{k'=1}^{K}\sum_{j\in\mathcal C_{k'}}l_{ij}\Gamma\bar{x}_{d}^{k'}\bigg]\\
&=&\sum_{i\in\mathcal
C_{k}}d_{i}(x^{i}-\bar{x}^{k}_{d})^{\top}\bigg[f_{k}(x^{i})
-f_{k}(\bar{x}^{k}_{d})
+c\sum_{k'=1}^{K}\sum_{j\in\mathcal{C}_{k'}}l_{ij}\Gamma
(x^{j}-{\bar{x}}^{k'}_{d})\bigg]
\end{eqnarray*}

From the decreasing condition (\ref{decreasing}):
\begin{eqnarray*}
(w-v)^{\top}[f_{k}(w)-f_{k}(v)-\alpha\Gamma(w-v)]\le-\delta(w-v)^{\top}(w-v),
\end{eqnarray*}
we have
\begin{eqnarray*}
\dot{V}_{k} &\le&-\delta\sum_{i\in \mathcal{C}_k}d_{i}
(x^{i}-\bar{x}^{k}_{d})^{\top}(x^{i}-\bar{x}^{k}_{d})\\
&&+\sum_{i\in\mathcal C_{k}}d_{i}(x^{i}-\bar{x}^{k}_{d})^{\top}
\bigg[c\sum_{k'=1}^{K}\sum_{j\in\mathcal{C}_{k'}}l_{ij}\Gamma
(x^{j}-{\bar{x}}^{k'}_{d})+\alpha\Gamma(x^{i}-\bar{x}^{k}_{d})\bigg].
\end{eqnarray*}
Thus,
\begin{eqnarray*}
\dot{V}&\le&-\delta\sum_{k=1}^{K}\sum_{i\in\mathcal
C_{k}}d_{i}(x^{i}-\bar{x}^{k}_{d})^{\top}(x^{i}-\bar{x}^{k}_{d})\\
&&+\sum_{k=1}^{K}\sum_{i\in\mathcal
C_{k}}d_{i}(x^{i}-\bar{x}^{k}_{d})^{\top}
\bigg[c\sum_{k'=1}^{K}\sum_{j\in\mathcal{C}_{k'}}l_{ij}\Gamma
(x^{j}-{\bar{x}}^{k'}_{d})+\alpha\Gamma(x^{i}-\bar{x}^{k}_{d})\bigg]\\
&=&-\delta\sum_{k=1}^{K}\sum_{i\in\mathcal
C_{k}}d_{i}(x^{i}-\bar{x}^{k}_{d})^{\top}(x^{i}-\bar{x}^{k}_{d})+(x-\bar{x}_{d})^{\top}\Big\{\big[D(cL+\alpha
I_{m})\big]^{s}\otimes \Gamma\Big\} (x-\bar{x}_{d})
\end{eqnarray*}
where $\otimes$ denotes the Kronecker product and $D={\rm
diag}[d_{1},\cdots,d_{m}]$.

It is clear that $[D(cL+\alpha I_{m})]^{s}\bigg|_{\mathcal
T_{\mathcal C}^{d}(1)}\le 0$ implies $\big\{[D(cL+\alpha
I_{m})]^{s}\otimes I_{n}\big\}\bigg|_{\mathcal T_{\mathcal
C}^{d}(n)}\le 0$. Decompose the positive definite matrix $\Gamma$ as
$\Gamma=C^{\top}C$ for some matrix $C$ and let
$y=[{y^{1}}^{\top},\cdots,{y^{m}}^{\top}]^{\top}$ with
$y^{i}=C(x^{i}-\bar{x}_{d}^{k})$ for all $i\in\mathcal C_{k}$, i.e.,
$y=(I_{m}\otimes C)(x-\bar{x}_{d})$. By Lemma 1, it is easy to see
that $ \sum_{i\in\mathcal C_{k}}d_{i}y^{i}=\sum_{i\in\mathcal
C_{k}}d_{i}C(x^{i}-\bar{x}_{d}^{k})=0$. This implies that
$y\in\mathcal T_{\mathcal C}^{d}(n)$. Therefore,
\begin{eqnarray}
&&(x-\bar{x}_{d})^{\top}\Big\{\big[D(cL+\alpha
I_{m})\big]^{s}\otimes
\Gamma\Big\}(x-\bar{x}_{d})\nonumber\\
&&=(x-\bar{x}_{d})^{\top}\Big(I_{m}\otimes
C^{\top}\Big)\Big\{\big[D(cL+\alpha I_{m})\big]^{s}\otimes
I_{n}\Big\}\Big(I_{m}\otimes C\Big)(x-\bar{x}_{d})\nonumber\\
&&=y^{\top}\Big\{\big[D(cL+\alpha I_{m})\big]^{s}\otimes
I_{n}\Big\}y\le 0.\label{form1}
\end{eqnarray}
Hence, we have
$$\dot{V}\le-\delta (x-\bar{x}_{d})^{\top}(D\otimes
I_{n})(x-\bar{x}_{d})=-2\delta\times V.$$ This implies that
$V(t)\le\exp(-2\delta~t)V(0)$. Therefore, $\lim_{t\to\infty}V(t)=0$.
Namely. $\lim_{t\to\infty}[x(t)-\bar{x}_{d}(t)]=0$ holds. In other
words, $\lim_{t\to\infty}[x^{i}-\bar{x}_{d}^{k}]=0$ for each
$i\in\mathcal C_{k}$ and $k=\onetoK$. According to the assumption
that $f_{k}(\cdot)$ are so different that if  cluster
synchronization is realized, the clusters are also different, we are
safe to say that the coupled system (\ref{Eq.4}) can clustering
synchronize.  $\square$

If each uncoupled system $\dot{x}^{i}=f_{k}(x^{i})$ is unstable, in
particular, chaotic, $\alpha$ must be positive in the inequality
(\ref{decreasing}). It is natural to raise the question: can we find
some positive diagonal matrix $D$ such that (\ref{state1}) satisfies
with sufficiently large $c$ and some certain $\alpha>0$. In other
words, for the coupled system (\ref{Eq.3}), what kind of unweighted
graph topology $\mathcal G$ satisfying the common inter-cluster
condition (\ref{inter-cluster}) can be a {\em chaos cluster
synchronizer} with respect to the clustering $\mathcal C$. It can be
seen that if the restriction of $(DL+L^{\top}D)$ to the transverse
subspace $\mathcal T_{\mathcal C}^{d}(1)$ is negative, i.e.,
\begin{eqnarray}
(DL+L^{\top}D)|_{\mathcal T_{\mathcal C}^{d}(1)}<0\label{Ineq1}
\end{eqnarray}
holds, then Ineq. (\ref{state1}) holds for sufficiently large $c$.

With these observations, we have

\begin{theorem}\label{state2}
Suppose that the common inter-cluster coupling condition
(\ref{inter-cluster}) holds for the coupled system (\ref{Eq.4}), and
$\alpha>0$. There exist a positive diagonal matrix $D$ and a
sufficiently large constant $c$ such that Ineq. (\ref{state1}) holds
if and only if all vertices in the same cluster belong to the same
connected component \footnote{Connected component in a direct graph
$\mathcal G$ is a maximal vertex set of which each vertices can
access all others.} in the graph $\mathcal G$.
\end{theorem}
{\em Proof}.  We prove the sufficiency for connected graph and
unconnected graph separated.

{\em Case 1}: The graph $\mathcal G$ is connected. Then, $L$ is
irreducible. Perron-Frobenius theorem (see Ref. \cite{Hor} for more
details) tells that the left eigenvector
$[\xi_{1},\cdots,\xi_{m}]^{T}$ of $L$ associated with the eigenvalue
$0$ has all components $\xi_{i}>0$, $i=1,\cdots,m$. In this case, we
pick $d_{i}=\xi_{i}$, $i=1,\cdots,m$. And its symmetric part
$[DL]^{s}=(DL+L^{\top}D)/2$ has all row sums zero and irreducible
with $\lambda_{1}([DL]^{s})=0$ associated with the eigenvector
$e=[1,\cdots,1]^{\top}$ and $\lambda_{2}([DL]^{s})<0$. Therefore,
$u^{\top}(DL)u\leq \lambda_{2}(DL)^{s}u^{\top}u<0$ for any $u\ne 0$
satisfying $u^{\top}e=0$.

Now, for any $u=[u_{1},\cdots,u_{m}]^{\top}\in\mathbb R^{m}$ with
$u^{\top}d=0$, define $\tilde{u}=[\bar{u},\cdots,\bar{u}]^{\top}$,
where $\bar{u}=\frac{1}{m}\sum_{i=1}^{m}u_{i}$. It is clear that
$DL\tilde{u}=0$ and $\tilde{u}^{\top}DL=0$ and
$(u-\tilde{u})^{\top}e=0$. Therefore,
\begin{eqnarray*}
u^{\top}(DL+L^{\top}D)u=(u-\tilde{u})^{\top}(DL+L^{\top}D)(u-\tilde{u})<0
\end{eqnarray*}
since  both hold. This implies Ineq. (\ref{Ineq1}) holds.

{\em Case 2}: The graph $\mathcal G$ is disconnected. In this case,
we can divide the bigraph $\mathcal G$ into several connected
components. If all vertices belongs to the same cluster are in the
same connected component, then by the same discussion done in case
1, we conclude that Ineq. (\ref{Ineq1}) holds for some positive
definite diagonal matrix $D$.

Necessity. We prove the necessity by reduction to absurdity.
Considering an arbitrary disconnected graph $\mathcal G$, without
loss of generality, supposing that $L$ has form:
\begin{eqnarray*}
L=\left[\begin{array}{ll}L_{1}&0\\
0&L_{2}\end{array}\right],
\end{eqnarray*}
and letting $\mathcal V_{1}$ and $\mathcal V_{2}$ correspond to the
sub-matrices $L_{1}$ and $L_{2}$ respectively, we assume that there
exists a cluster $\mathcal C_{1}$ satisfying $\mathcal
C_{1}\bigcap\mathcal V_{i}\ne\emptyset$ for all $i=1,2$. That is,
there exists at least a pair of vertices in the cluster $\mathcal
C_{1}$ which can not access each other. For each
$d=[d_{1},\cdots,d_{m}]^{\top}$ with $d_{i}>0$ for all $i=\onetom$,
letting $D={\rm diag}[d_{1},\cdots,d_{m}]$, we can find a nonzero
vector $u\in\mathcal T_{\mathcal C}^{d}(1)$ such that $u^{\top}DL
u=0$ (see the Appendix for the details). This implies that
inequality (\ref{Ineq1}) does not hold. So, inequality
(\ref{state1}) can not hold for any positive $\alpha$. $\square$

In the case that the clustering synchronized trajectories are
chaotic, with $\alpha>0$, Theorem \ref{state2} tells us that chaos
cluster synchronization can be achieved (for sufficiently large
coupling strength) if and only if all vertices in the same cluster
belongs to the same connected component in the graph $\mathcal G$.

In summary, the following two conditions play the key role in
cluster synchronization.
\begin{enumerate}
\item Common inter-cluster edges for each vertex in the same cluster;

\item Communicability for each pair of vertices in the same cluster.
\end{enumerate}
The first condition guarantees that the clustering synchronization
manifold is invariant through the dynamical system with properly
picked weights and the second guarantees that chaos clustering
synchronization can be reached with a sufficiently large coupling
strength.

\subsection{Schemes to clustering synchronize}
The theoretical results in previous section indicate that the
communication among vertices in the same cluster is important for
chaos cluster synchronization. A cluster is said to be {\em
communicable} if every vertex in this cluster can connect any other
vertex by paths in the global graph. These paths between vertices
are composed of edges, which can be either of inter-cluster or
intra-cluster. Refs. \cite{Jalan} showed that this classification of
paths distinguishes the formation of clusters. A self-organized
cluster implies that the intra-cluster edges are dominant for the
communications between vertices in this cluster. And, a driven
cluster is that the inter-cluster edges are dominant for the
communications between vertices in this cluster. There are various
of ways to describe ``domination''. In the following, we consider
the unweighted graph topology and investigate the two clustering
schemes via the results above.

\begin{figure}
\includegraphics[width=15cm]{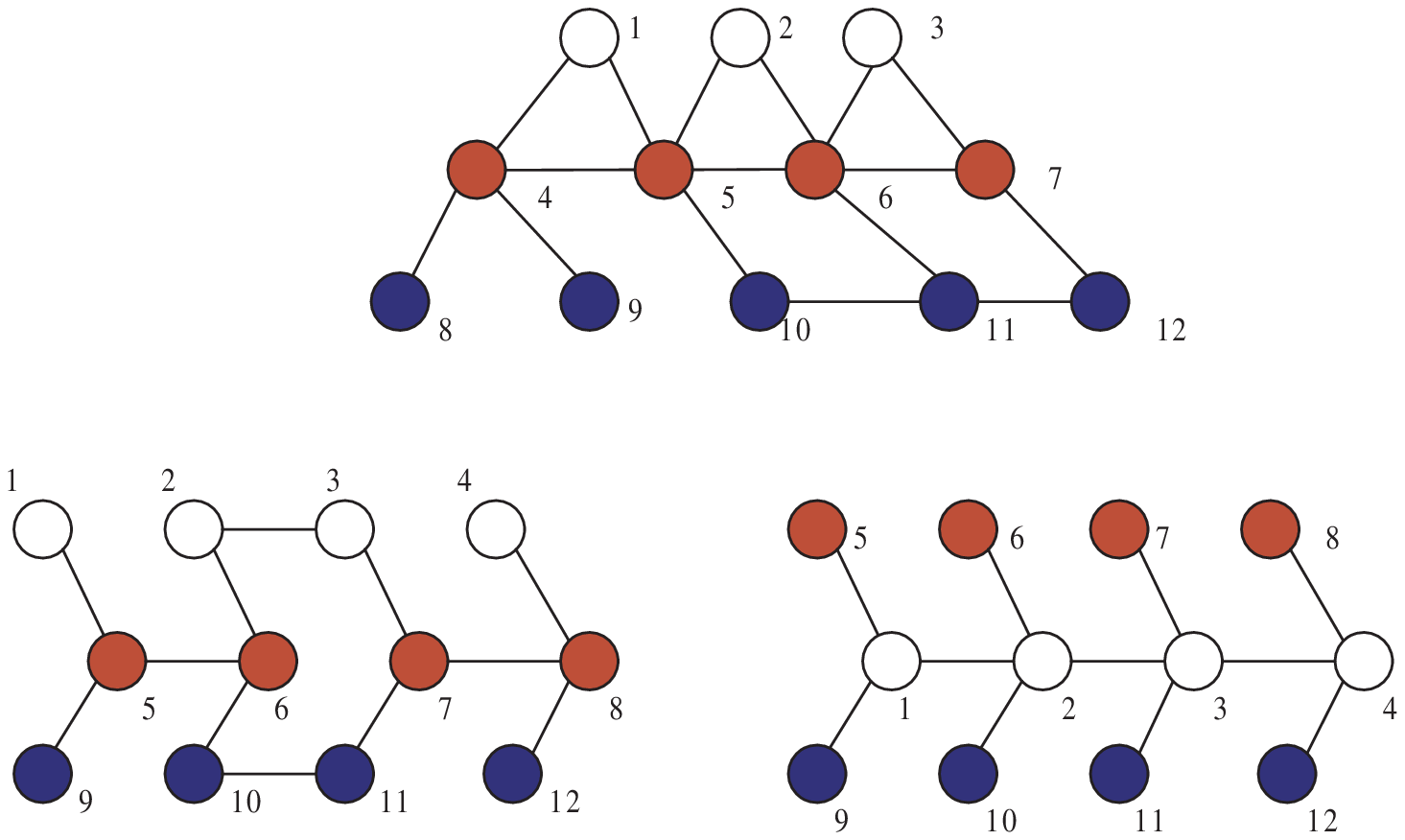}
\caption{\label{graphs}Graphs of examples. In the graph 1, the white
cluster (vertex set $1-3$) is driven since they has no intra-cluster
edges, the red cluster (vertex set $4-7$) is mixed since each pari
of vertices can access each other via only inter- or intra- edges,
and the blue cluster is driven since each pair of vertices can
access each other via only the inter-cluster edges but can not
communicate only via intra-cluster edges. In the  graph 2, each
cluster of the white and blue clusters (vertex sets $1-4$ and
$9-12$) is driven since each pair of vertices can access each other
only via inter-cluster edges but only has a single intra-cluster
edges. However, the read cluster (vertex set $5-8$) is recognized as
a hybrid cluster since the sets of inter- or intra-cluster edges are
both necessary for communication between each pair of vertices. In
the graph 3, the red and blue clusters (vertex sets $5-8$ and
$9-12$) are all driven since they do not have intra-cluster edges
and the white cluster (vertex set $1-4$) is an example of
self-organization since each pair of vertices can communicate via
only the intra-cluster edges but can not if removing the
intra-cluster edges.}
\end{figure}

We firstly describe {\em self-organization} and {\em driving} as two
schemes for cluster synchronization. Self-organization represents
the scheme that the set of intra-cluster edges are irremovable for
the communication between each pair of vertices in the same cluster
and driving represents the scheme that the set of inter-cluster
edges are irremovable for the communication between vertices in the
same cluster. Thus, we propose the following classification of
clusters.
\begin{enumerate}

\item {\em Self-organized cluster}: the subgraph of the cluster is
connected but if removing the intra-cluster links of the cluster,
there exist at least one pair of vertices such that no paths in the
remaining graph can connect them;

\item {\em Driven cluster}: the subgraph of the cluster is
disconnected but even if removing all intra-cluster links of the
cluster, each pair of vertices in the cluster can reach each other
by paths in the remaining graph;

\item {\em Mixed cluster}: the subgraph of the cluster is connected
and even if removing all intra-cluster links of the cluster, each
pair of vertices in the cluster can reach each other by paths in the
remaining graph;

\item {\em Hybrid cluster}: the subgraph of the cluster is
disconnected and if removing the intra-cluster links of the cluster,
there exist at least one pair of vertices such that no paths in the
remaining graph can link them.

\end{enumerate}
Table \ref{class} describes the characteristics of each cluster
class.
\begin{table}

\renewcommand{\arraystretch}{1.3}
\caption{Communicability of clusters under edge-removing operations}
\label{class}

\centering
\begin{tabular}{|c|c|c|}
\hline &Remove the intra-cluster edges & Remove the inter-cluster
edges \\
\hline
Self-organized cluster & no&yes\\
\hline
Driven cluster  & yes&no\\
\hline
Mixed cluster&yes&yes\\
\hline Hybrid cluster&no&no\\
 \hline
\end{tabular}
\end{table}
Fig.\ref{graphs} shows examples of these four kinds of clusters,
which will be used in later numerical illustrations. With this
cluster classification, we conclude that any mixed or self-organized
cluster can not access another hybrid or self-organized cluster.
Table \ref{inter} shows all possibilities of accessibility among all
kinds of clusters in a connected graph. Moreover, it should be
noticed that the cluster in the networks as illustrated in
Fig.\ref{graphs} may not be connected via the subgraph topologies.
For example, the white and blue clusters in graph 1, the red and
blue clusters in the graph 3, as well as all clusters in graph 2,
are not connected by inter-cluster subgraph topologies. Certainly,
the vertices in the same cluster are connected via inter- and
intra-cluster edges. That is, we can realize cluster synchronization
in non-clustered networks.
\begin{table}

\renewcommand{\arraystretch}{1.3}
\caption{Possibility of coexistence for two kinds of clusters in
connected graph} \label{inter}

\centering
\begin{tabular}{|c|c|c|c|c|}
\hline &Self-organized & Driven & Hybrid & Mixed \\
\hline
Self-organized & $\times$& $\surd$ & $\times$& $\times$\\
\hline
Driven & $\surd$&$\surd$&$\surd$&$\surd$\\
\hline
Hybrid&$\times$&$\surd$&$\surd$&$\times$\\
\hline Mixed&$\times$&$\surd$&$\times$&$\surd$\\
 \hline
\end{tabular}
\end{table}

\subsection{Examples}
In this part, we propose several numerical examples to illustrate
the theoretical results. In this example, $K=3$. The three graph
topologies are shown in Fig. \ref{graphs}. The coupled system is
\begin{eqnarray}
\dot{x}^{i}=f_{k}(x^{i})
+c\bigg[\sum_{\mathcal{N}_{k'}(i)\ne\emptyset}\frac{1}{d_{i,k'}}
\sum_{j\in\mathcal N_{k'}(i)}\Gamma(x^{j}-x^{i})\bigg],~i\in\mathcal
C_{k},~k=1,2,3,\label{simu1}
\end{eqnarray}
where $\Gamma={\rm diag}[1,1,0]$ and $f_{k}(\cdot)$ are
non-identical Lorenz systems:
\begin{eqnarray}
f_{k}(u)=\left\{\begin{array}{l}
10(u_{2}-u_{1})\\
\displaystyle\frac{8}{3}u_{1}-u_{2}-u_{1}u_{3}\\
u_{1}u_{2}-b_{k}u_{3},
\end{array}\right.\label{Lorenz}
\end{eqnarray}
where the parameter $b_{1}=28$ for the white cluster, $b_{2}=38$ for
the red cluster, and $b_{3}=58$ for the blue cluster.

\begin{figure}
\includegraphics[width=15cm]{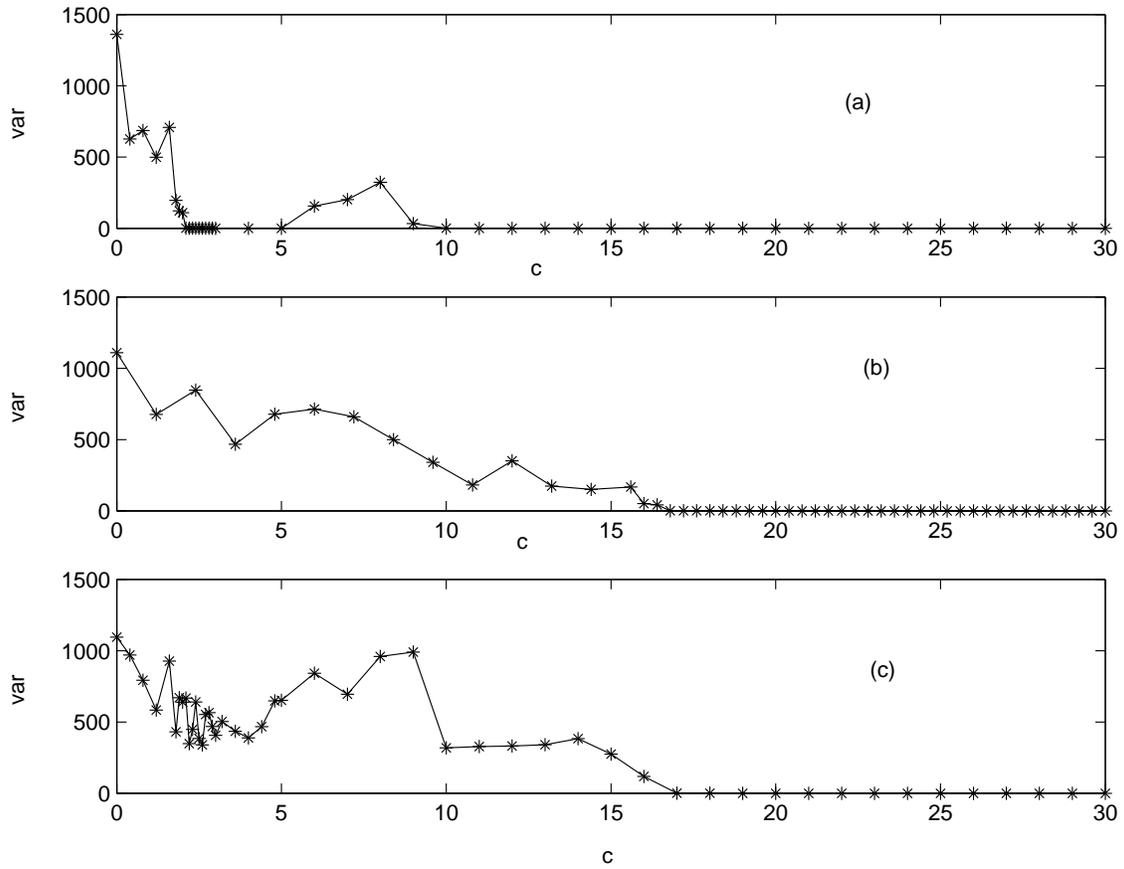}
\caption{\label{var} $var$ with respect to $c$ : (a) for the graph
1; (b) for the graph 2; (c) for the graph 3, respectively.}
\end{figure}

As shown in Ref. \cite{XC}, the boundedness of the trajectories of
an array of coupled Lorenz systems can be ensured. Therefore, the
decreasing condition (\ref{decreasing}) is satisfied for a
sufficiently large $\alpha$. We use the following quantity to
measure the variation for vertices in the same cluster:
\begin{eqnarray*}
{\rm var}=\bigg\langle\sum_{k=1}^{K}\frac{1}{\# \mathcal
C_{k}-1}\sum_{i\in\mathcal
C_{k}}[x^{i}-\bar{x}_{k}]^{\top}[x^{i}-\bar{x}_{k}]\bigg\rangle
\end{eqnarray*}
where $\bar{x}_{k}=\frac{1}{\#\mathcal C_{k}}\sum_{i\in\mathcal
C_{k}}x^{i}$, $\langle\cdot\rangle$ denotes the time average. The
ordinary differential equations (\ref{simu1}) are solved by the
Runge-Kutta fourth-order formula with a step length 0.01. The time
average interval is $[50,100]$ \cite{interval}. Fig. \ref{var}
indicates that for either graph 1, graph 2, or graph 3, the coupled
system (\ref{simu1}) clustering synchronizes respectively, if the
coupling strength is larger than certain threshold value.  Instead,
for the graph 1, despite the coupled system can synchronize if $c$
is greater than some value (around 10), it can also synchronize if
$c\in[2.2,5]$). It is not very surprising. Previous theoretical
results only give sufficient condition that the coupled system can
clustering synchronize if the coupling strength $c$ is large enough.
It does not exclude the case that the coupled system can still
clustering synchronize even if the coupling strength $c$ is small.

The following quantity is used to measure the deviation between
clusters:
\begin{eqnarray*}
{\rm dis}(t)= \min_{i\ne
j}[\bar{x}_{i}(t)-\bar{x}_{j}(t)]^{\top}[\bar{x}_{i}(t)-\bar{x}_{j}(t)]
\end{eqnarray*}

\begin{figure}
\includegraphics[width=15cm]{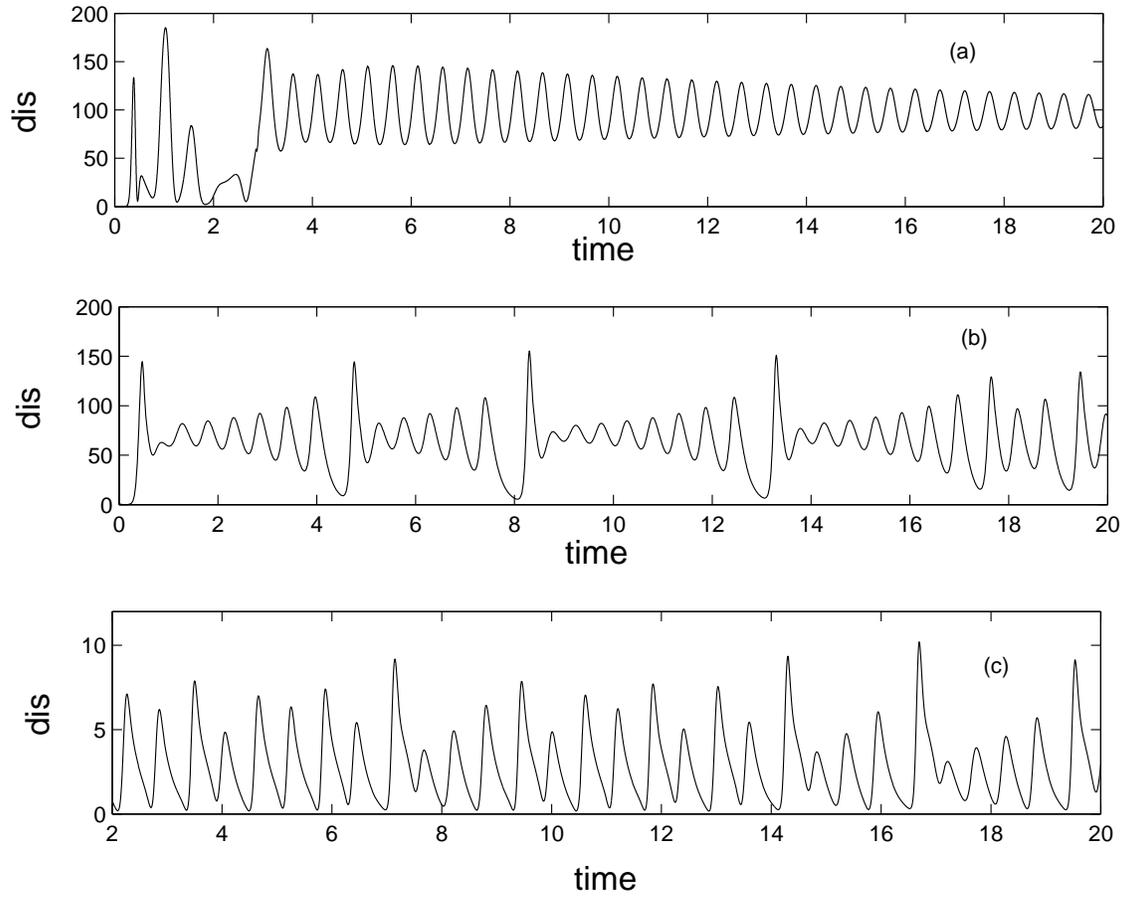}
\caption{\label{dis} Dynamics of $dis(t)$ through Eq. (\ref{Eq.3}):
(a) for the graph 1 with $c=12$; (b) for the graph 2 with $c=25$;
(c) for the graph 3 with $c=20$, respectively.}
\end{figure}

\begin{figure}
\includegraphics[width=15cm]{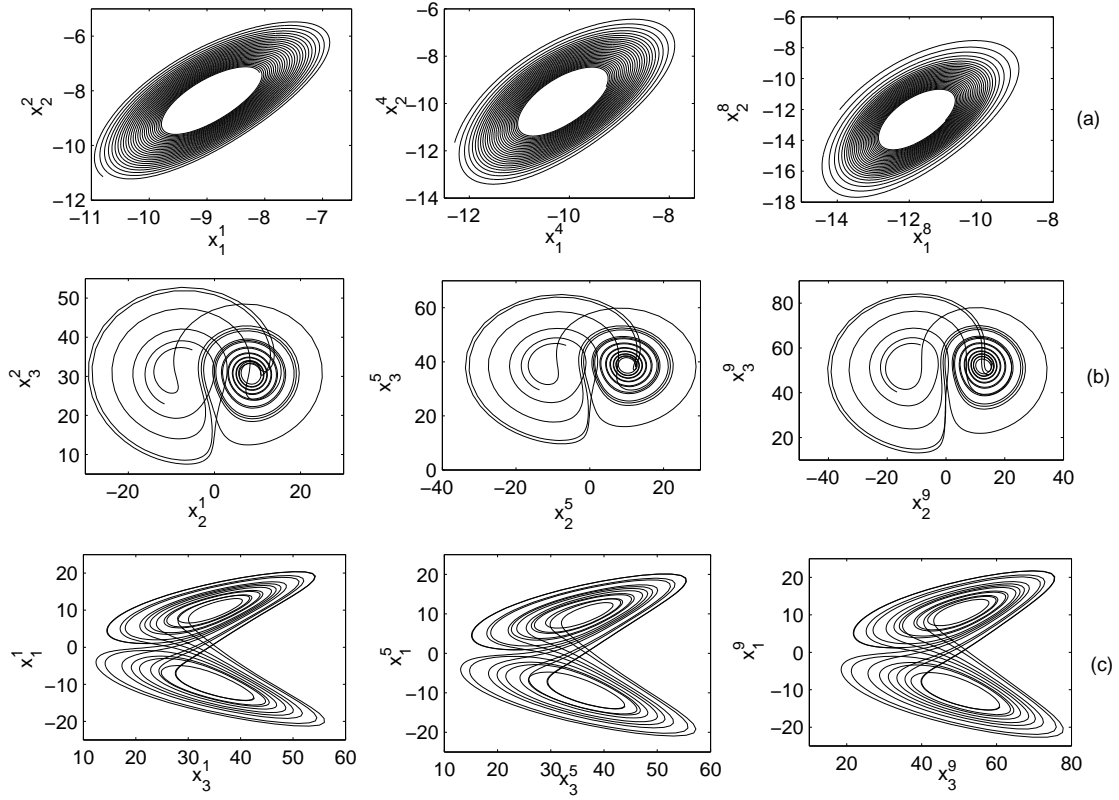}
\caption{\label{phase} Phase dynamics for each cluster through
equality (\ref{Eq.3}): (a) $x_{1}-x_{2}$ phase dynamics for each
cluster in the graph 1; (b) $x_{2}-x_{3}$ phase dynamics for each
cluster in the graph 2; (2) $x_{3}-x_{1}$ phase dynamics for each
cluster in the graph 3, respectively. }
\end{figure}

Fig.\ref{dis} shows that the deviation between clusters is apparent,
even "${\rm var}\simeq 0$". In Fig.\ref{phase}, the dynamical
behaviors for all clusters in certain phase plane are given.
Although the attractor for each cluster seems to have similar
structure and shape, the positions at same time are still different.
It is clear that the difference is caused by the different choice of
parameters for different clusters. This illustrates that the cluster
synchronization is actually realized.

\section{Adaptive feedback cluster synchronization algorithm}

For certain network topology which has weak cluster
synchronizability, i.e., the threshold to ensure clustering
synchronization is relatively large, which is further studied in
Sec. IV.A. It is natural to raise the following question: How to
achieve cluster synchronization for networks no matter whether they
have "good" topology or not. One approach proposed recently is
adding weights to vertices and edges. Refs \cite{Weight} showed
evidences that certain weighting procedures can actually enhance
complete synchronization. On the other hand, adaptive algorithm has
emerged as an efficient means of weighting to actually enhance
complete synchronizability \cite{Adaptive}.

In this
section, we consider the coupled system 
\begin{eqnarray}
\dot{x}^{i}=f_{k}(x^{i})+\sum_{j=1}^{m}a_{ij}w_{ij}\Gamma(
x^{j}-x^{i}),~i\in\mathcal C_{k},~k=\onetoK.\label{Eq.21}
\end{eqnarray}
and propose an adaptive feedback algorithm to achieve cluster
synchronization for a prescribed graph.

Suppose that the common inter-cluster and communicability conditions
are satisfied. Without loss of generality, we suppose that the graph
$G$ is undirected and connected. Consider the coupled system
(\ref{Eq.2}) with Laplacian $L$ defined as in Eqs.(\ref{Eq.4}) and
$d^{\top}=[d_{1},\cdots,d_{m}]$ is the left eigenvector of $L$
associated with the eigenvalue $0$.

Now, we propose the following adaptive cluster synchronization
algorithm
\begin{eqnarray}
\left\{\begin{array}{ll}
\dot{x}^{i}(t)=f_{k}(x^{i}(t))+\sum_{j=1}^{m}a_{ij}w_{ij}(t)\Gamma[
x^{j}(t)-x^{i}(t)],~i\in\mathcal C_{k},~k=\onetoK,\\
\dot{w}_{ij}(t)=\rho_{ij}
d_{i}[x^{i}(t)-\bar{x}^{k}_{d}(t)]^{\top}\Gamma[x^{i}(t)-x^{j}(t)],\\~{\rm
for~each} ~e_{ij}\in\mathcal E~{\rm and~}i\in\mathcal C_{k},
~k=\onetoK\end{array}\right.\label{adaptive}
\end{eqnarray}
with $\rho_{ij}>0$ are constants.
\begin{theorem}
Suppose that  the graph $\mathcal G$ is connected, all the
assumptions of Theorem \ref{st1} hold, the system (\ref{adaptive})
is essential bounded. The system (\ref{adaptive}) clustering
synchronizes for any initial data.
\end{theorem}
{\em Proof}. First of all, pick $l_{ij}$ as defined in Eqs.
(\ref{Eq.4}) and a sufficiently large $c$.Since $\mathcal G$ is
connected, Theorem 2 tells
\begin{eqnarray}
\bigg[D(cL+\alpha I_{m})\bigg]^{s}\bigg |_{\mathcal T_{\mathcal
C}^{d}(1)}< 0.\label{Ineq2}
\end{eqnarray}
Define the following candidate Lyapunov function
\begin{eqnarray*}
Q_{k}(x,W)=\sum_{i\in\mathcal
C_{k}}\bigg[\frac{d_{i}}{2}(x^{i}-\bar{x}^{k}_{d})^{\top}(x^{i}-\bar{x}^{k}_{d})
+\frac{1}{2\rho_{ij}}a_{ij}(w_{ij}-cl_{ij})^{2}\bigg],~
Q(x,W)=\sum_{k=1}^{K}Q_{k}.
\end{eqnarray*}
Differentiating $Q_{k}$, we have

\begin{eqnarray*}
\dot{Q}_{k} &=&\sum_{i\in\mathcal
C_{k}}d_{i}(x^{i}-\bar{x}^{k}_{d})^{\top}\bigg\{f_{k}(x^{i})
+\sum_{j=1}^{m}a_{ij}w_{ij}\Gamma(x^{j}-\bar{x}^{j})\bigg\}\\
&&+\sum_{i\in\mathcal C_{k}}\sum_{k'=1}^{K}\sum_{j\in\mathcal
N_{k'}(i)}a_{ij}(w_{ij}-cl_{ij})
d_{i}(x^{i}-\bar{x}_{d}^{k})^{\top}\Gamma(x^{i}-x^{j})\\
&=&\sum_{i\in\mathcal
C_{k}}d_{i}(x^{i}-\bar{x}^{k}_{d})^{\top}\bigg\{f_{k}(x^{i})
+c\sum_{j=1}^{m}l_{ij}\Gamma(x^{j}-x^{i})
-\dot{\bar{x}}_{d}^{k}\bigg\}.
\end{eqnarray*}
Similar to the proof of Theorem 1, we have
\begin{eqnarray*}
\sum_{i\in\mathcal C_{i}}\dot{Q}_{i}=\sum_{i\in\mathcal
C_{i}}d_{i}(x^{i}-\bar{x}^{k}_{d})^{\top}\bigg\{f_{k}(x^{i})-f(\bar{x}_{d}^{k})
+c\sum_{j=1}^{m}l_{ij}\Gamma(x^{j}-\bar{x}_{d}^{j})\bigg\}
\end{eqnarray*}
and
\begin{eqnarray*}
\dot{Q}&=&\sum_{k=1}^{K}\dot{Q}_{k}\le-\delta\sum_{k=1}^{K}\sum_{i\in\mathcal
C_{k}}d_{i}(x^{i}-\bar{x}^{k}_{d})^{\top}(x_{i}-\bar{x}^{k}_{d})\\&&+
\sum_{k=1}^{K}\sum_{i\in\mathcal
C_{k}}d_{i}(x^{i}-\bar{x}^{k}_{d})^{\top}
\bigg[\alpha\Gamma(x^{i}-\bar{x}^{k}_{d})
+c\sum_{j=1}^{m}l_{ij}\Gamma(x^{j}-\bar{x}_{d}^{j})\bigg]\\
&=&-\delta(x-\bar{x}_{d})^{\top}\big(D\otimes
I\big)(x-\bar{x}_{d})+(x-\bar{x}_{d})^{\top}\bigg\{[D(cL+\alpha
I_{m})]^{s}\otimes\Gamma\bigg\}(x-\bar{x}_{d}).
\end{eqnarray*}
Ineq \ref{Ineq2} implies
$$\dot{Q}\le-\delta(x-\bar{x}_{d})^{\top}\big(D\otimes
I\big)(x-\bar{x}_{d})\le 0.$$ This implies
\begin{eqnarray}
\int_{0}^{t}\delta(x(s)-\bar{x}_{d}(s))^{\top}\big(D\otimes
I\big)(x(s)-\bar{x}_{d}(s))ds\le Q(0)-Q(t)\le Q(0)<\infty \label{ad}
\end{eqnarray}
From the assumption of the boundedness of Eq.(\ref{adaptive}), we
can conclude $\lim_{t\to\infty}[x(t)-\bar{x}_{d}(t)]=0$ due to the
fact that $x(t)$ is uniform continuous. This completes the proof.
$\square$

For the disconnected situation, we can split the graph into several
connected components and deal with each connected component by the
same means as above.

\begin{figure}
\includegraphics[width=15cm]{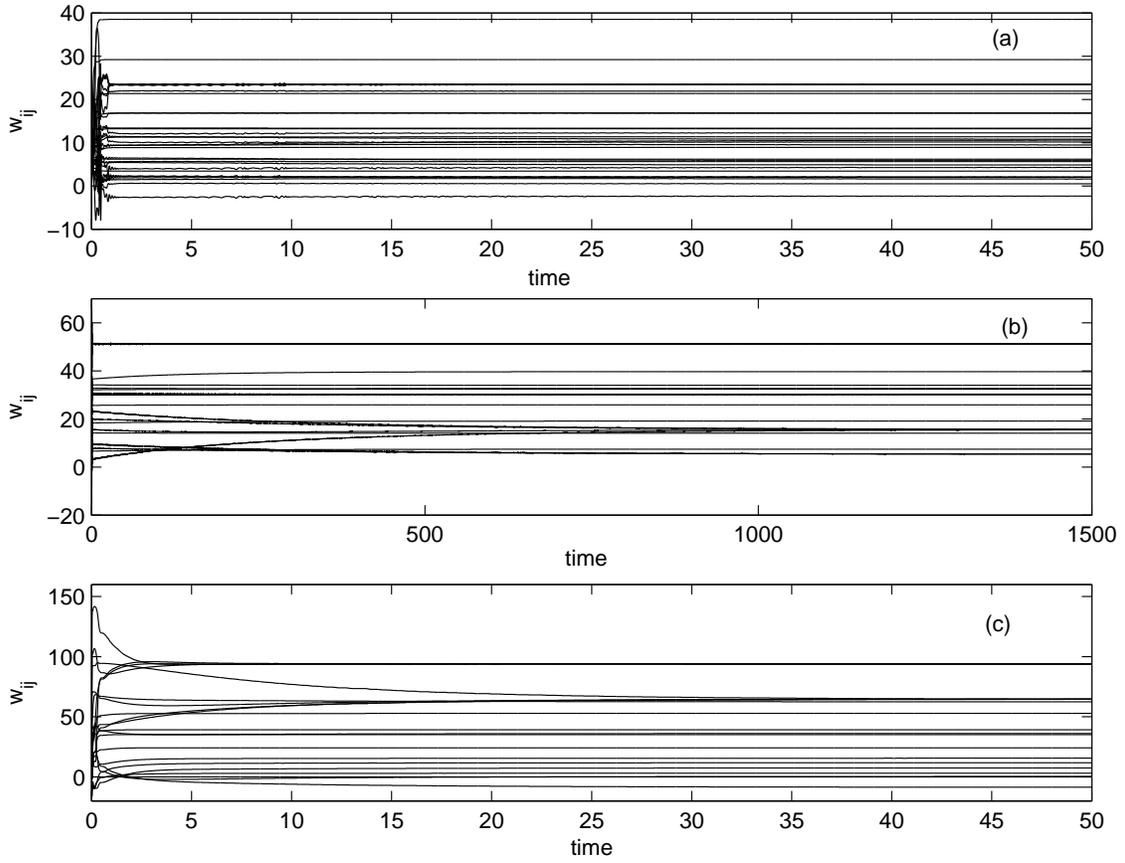}
\caption{\label{weights} Convergence dynamics of weights
$\{w_{ij},~(i,j)\in\mathcal E\}$ of edges through equality
(\ref{Eq.3}) with the adaptive algorithm (\ref{adaptive}): (a) for
the graph 1; (b) for the graph 2; (c) for the graph 3, respectively.
}
\end{figure}

The dynamics of the weights $w_{ij}(t)$ is an interesting issue.
Even though it is illustrated in Fig. \ref{weights} that all weights
converge, to our best reasoning, we can only prove that all
intra-weights converge, i.e., vertices $i$ and $j$ belonging to the
same cluster $\mathcal C_{k}$. In fact, by (\ref{ad}), we have
\begin{eqnarray*}
\int_{0}^{\infty} [x^{i}(\tau)-\bar{x}^{k}_{d}(\tau)]^{\top}
[x^{i}(\tau)-\bar{x}^{k}_{d}(\tau)]d\tau<+\infty.
\end{eqnarray*}
Thus,
\begin{eqnarray*}
&&\int_{0}^{\infty}|\dot{w}_{ij}(\tau)|d\tau
=\rho_{ij}d_{i}\int_{0}^{\infty}\bigg|[x^{i}(\tau)-\bar{x}^{k}_{d}(\tau)]
^{\top}\Gamma [x^{i}(\tau)-x^{j}(\tau)]\bigg|d\tau\\
&\le&\int_{0}^{\infty}\rho_{ij}d_{i}\|\Gamma\|_{2}
\bigg\{\big|[x^{i}(\tau)-\bar{x}^{k}_{d}(\tau)]^{\top}
[x^{i}(\tau)-\bar{x}^{k}_{d}(\tau)]\big|
+\big|[x^{i}(\tau)-\bar{x}^{k}_{d}(\tau)]^{\top}
[x^{j}(\tau)-\bar{x}^{k}_{d}(\tau)]\big|\bigg\}d\tau\\
&\le&\rho_{ij}d_{i}\|\Gamma\|_{2}\bigg\{\frac{3}{2}\int_{0}^{\infty}
[x^{i}(\tau)-\bar{x}^{k}_{d}(\tau)]^{\top}
[x^{i}(\tau)-\bar{x}^{k}_{d}(\tau)]d\tau\\
&&+\frac{1}{2}\int_{0}^{\infty}
[x^{j}(\tau)-\bar{x}^{k}_{d}(\tau)]^{\top}
[x^{j}(\tau)-\bar{x}^{k}_{d}(\tau)]d\tau\bigg\}.
\end{eqnarray*}
Therefore, for any $\epsilon>0$, there exists $T>0$, such that for
any $t_{1}>T$, $t_{2}>T$, we have
\begin{eqnarray*}
|w_{ij}(t_{2})-w_{ij}(t_{1})|\le\int_{t_{1}}^{t_{2}}|\dot{w}_{ij}(\tau)|d\tau
<\epsilon
\end{eqnarray*}
By Cauchy convergence principle, $w_{ij}(t)$ converges to some final
weights $w_{ij}^{*}$ for $i\in\mathcal{C}_{k}$,
$j\in\mathcal{C}_{k}$ when $t\rightarrow\infty$.

On the other hand, to our best reasoning, we can not prove whether
or not the weights $w_{ij}(t)$ converges, if the vertex $i$ and $j$
belong to different clusters. If we assume the convergence of all
weights, according to the the LaSalle invariant principle, the final
weights should guarantee that the cluster synchronization manifold
is still invariant. That is to say, if difference trajectories:
$s^{k'}-s^{k}$ in Eqs. (\ref{cluster_syn}), are linearly
independent, the cluster the condition (6) still holds for the final
weights.

\begin{figure}
\includegraphics[width=15cm]{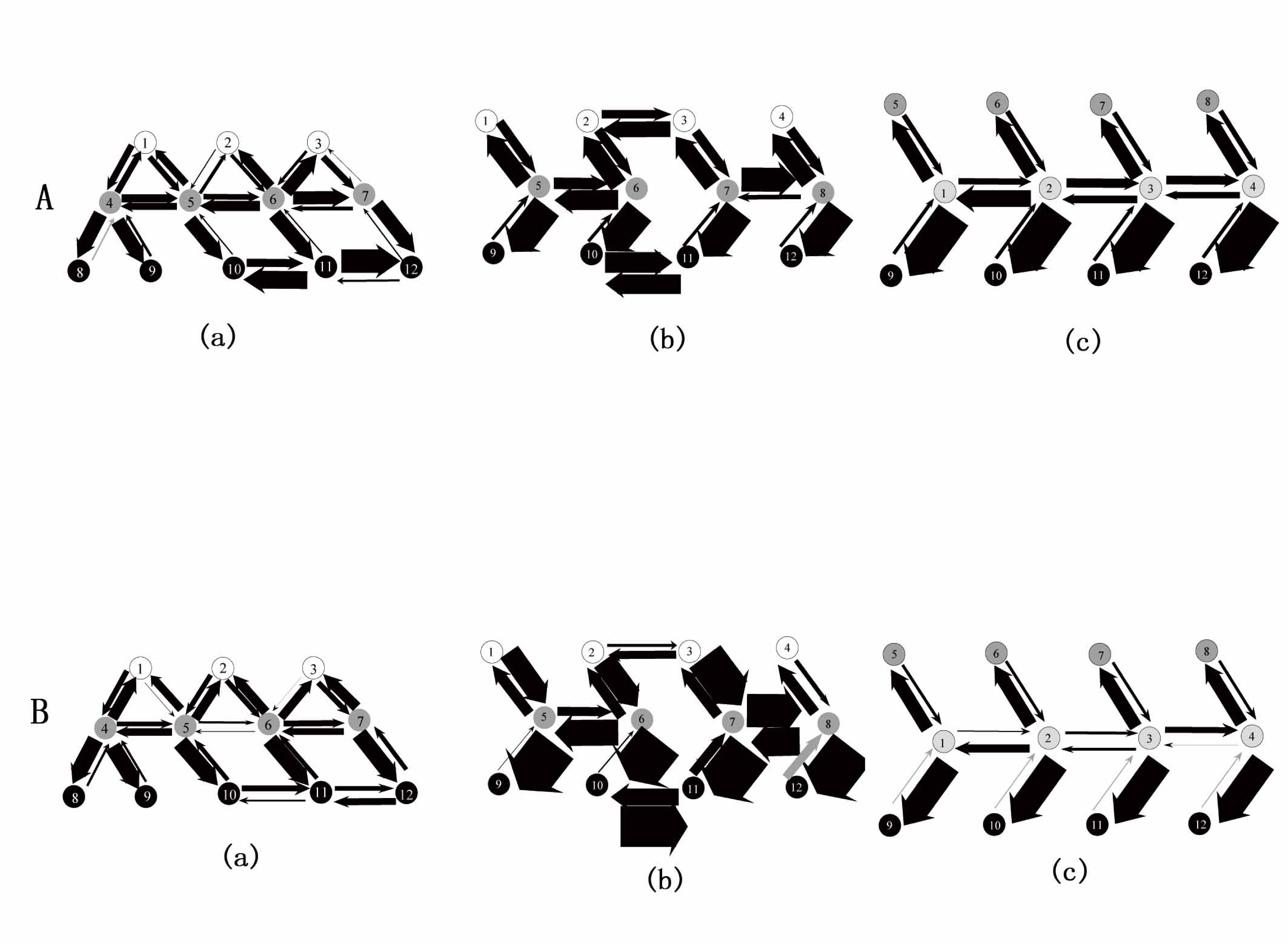}
\caption{\label{final_weights} Two sets of the final weighted
topologies of the three graphs in Fig. \ref{graphs} via employing
the adaptive algorithm (\ref{adaptive}) with two different sets of
initial data but the same parameters. Set A and B correspond two set
of initial values and (a)-(c) correspond the graph 1-3 in Fig.
\ref{graphs}. The color of the line represents the sign of the
weights (black for positive and gray for negative) and the width of
the line represents the scale of the weight in modulus.}
\end{figure}

Moreover, we have found out that the final weights in our example
sensitively depends on the initial values. Fig. \ref{final_weights}
gives two sets of weighted topologies of the three graphs as shown
in Fig. \ref{graphs} after employing the adaptive algorithm with two
different sets of initial values of $w_{ij}(0)$ and the same
parameters. One can see that the final weight can be quite different
for different initial values and even be negative. From this
observation, we argue that it may be the adaptive process not the
final weights counts to reach cluster synchronization. The further
investigation of the final weights is out of the scope of the
current paper.



\subsection{Examples}
\begin{figure}
\includegraphics[width=15cm]{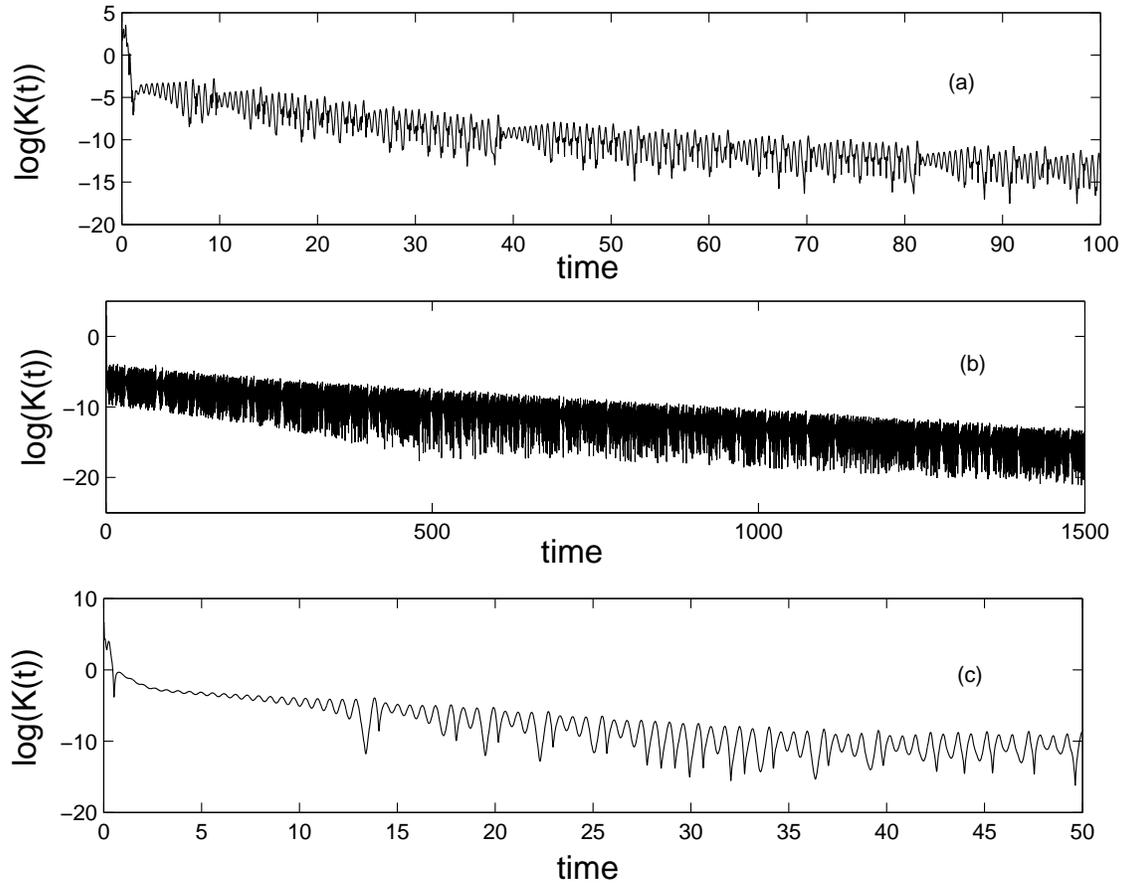}
\caption{\label{adp_K} Dynamics of the logarithm of $K(t)$ through
equality (\ref{Eq.3}) with the adaptive algorithm (\ref{adaptive}):
(a) for the graph 1; (b) for the graph 2; (c) for the graph 3,
respectively.}
\end{figure}

\begin{figure}
\includegraphics[width=15cm]{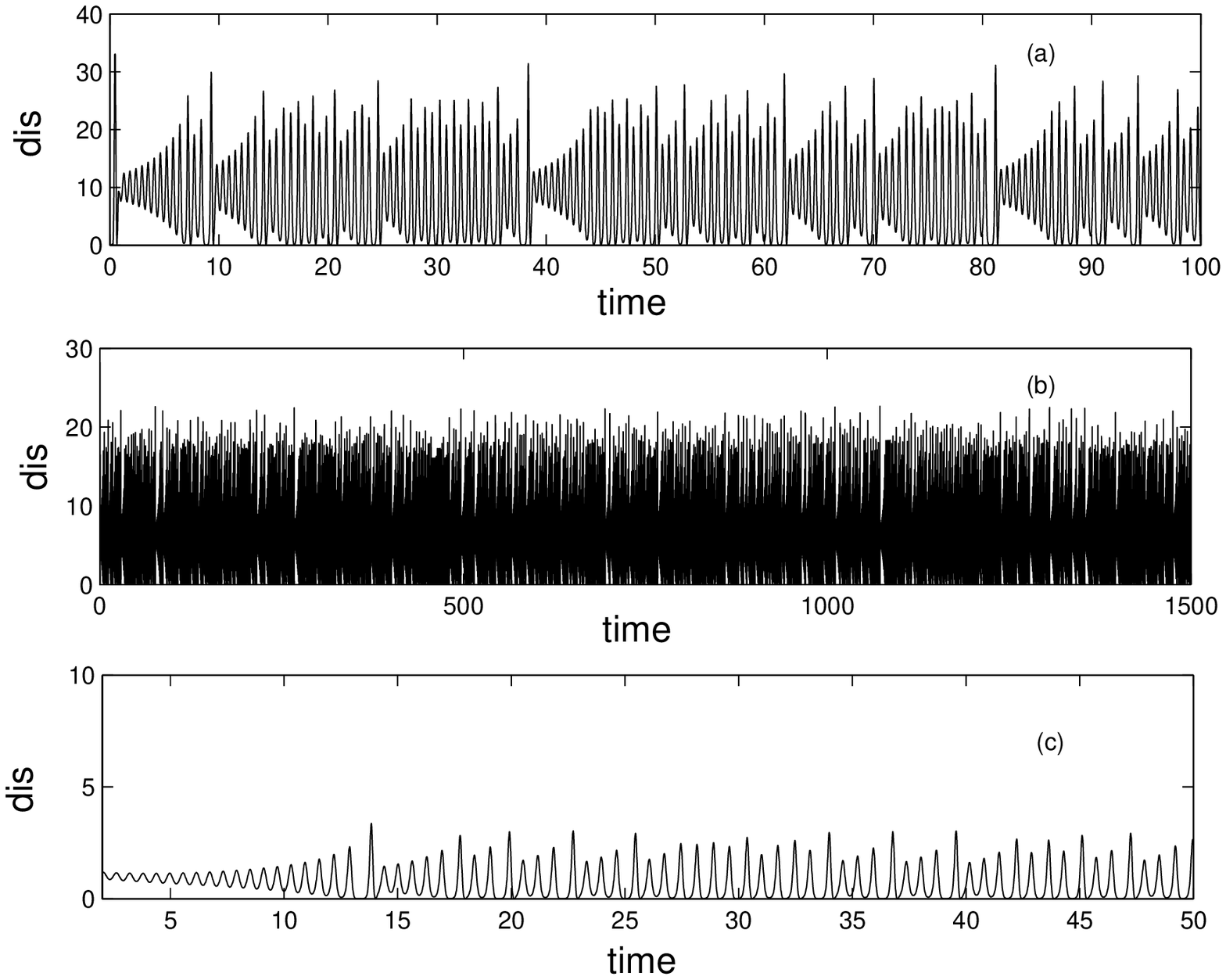}
\caption{\label{adp_dis} Dynamics of $dis(t)$ through equality
(\ref{Eq.3}) with the adaptive algorithm (\ref{adaptive}): (a) for
the graph 1; (b) for the graph 2; (c) for the graph 3,
respectively.}
\end{figure}

We still use the graphs 1-3 described in Fig.\ref{graphs} and the
Lorenz system (\ref{Lorenz}) as the uncoupled system to illustrate
the adaptive feedback algorithms. The ordinary differential
equations are solved by the the Runge-Kutta fourth-order formula
with a step length 0.01. The initial values of the states and the
weights are randomly picked in $[-3,3]$ and $[-5,5]$ respectively.
We use the following quantity to measure the state variance inside
each cluster with respect to time.
\begin{eqnarray*}
K(t)=\sum_{k=1}^{K}\frac{1}{\# \mathcal C_{k}-1}\sum_{i\in\mathcal
C_{k}}[x^{i}(t)-\bar{x}_{k}(t)]^{\top}[x^{i}(t)-\bar{x}_{k}(t)]
\end{eqnarray*}
Fig.\ref{adp_K} shows that the adaptive algorithm succeeds in
clustering synchronizing the network with respect to the pre-given
clusters. Figs.\ref{adp_dis} indicates that the differences between
clusters due to non-identical parameters $b_{k}$. As shown in
Fig.\ref{weights}, the weights converge but the limit values are not
always positive. This is not surprising. The right-hand side of the
algorithm (\ref{adaptive}) can be either positive or negative, which
causes some weights of edges to be negative. The situation of
negative weights is out of the scope of this paper.

\section{Discussions}
In this section, we make further discussions for some interesting
relating issues.

\subsection{Clustering synchronizability}
Synchronizability is used to measure the capability of of
synchronization for the graph. It can be described by the threshold
of the coupling strength to guarantee that the coupled system can
synchronize. For complete synchronization, it was formulated as a
function of the eigenvalues of symmetric Laplacian \cite{Loc_syn} or
certain Rayleigh quotient of asymmetric Laplacian \cite{Wu2}. How
the topology of the underlying graph affects synchronizability is an
important issue for the study of complex networks \cite{Boc1}. Here,
similarly, we are also interested in how to formulate and analyze
the cluster-synchronizability of a graph $\mathcal G$ and a
clustering $\mathcal C$.

Consider the model (\ref{Eq.4}) of coupled system. Theorem \ref{st1}
tells us that under the common inter-cluster condition, the cluster
synchronization condition (\ref{state1}) can be rewritten as
\begin{eqnarray}
c>\frac{\alpha}{\min_{u\in\mathcal T_{\mathcal C}^{d}(1),~u\ne
0}\frac{-u^{\top}(DL)^{s}u}{u^{\top}Du}}
\end{eqnarray}
for some positive definite diagonal $D$. Therefore, we take the
Rayleigh-Hitz quotient
\begin{eqnarray*}
CS_{\mathcal G,\mathcal C}=\max_{D\in\mathcal D}\min_{u\in\mathcal
T_{\mathcal C}^{d}(1),~u\ne 0}\frac{-u^{\top}(DL)^{s}u}{u^{\top}D
u},
\end{eqnarray*}
to measure the cluster synchronizability for the graph $\mathcal G$
and clustering $\mathcal C$, where $\mathcal D$ denotes the set of
positive definite diagonal matrices of dimension $m$. It can be seen
that the larger $CS_{\mathcal G,\mathcal C}$ is, the smaller the
coupling strength $c$ can be such that the coupled system
(\ref{Eq.4}) clusteringly synchronize. In particular, if $L$ is
symmetric, then $CS_{\mathcal G,\mathcal C}$ is just the maximum
eigenvalue of $-L$ in the transverse space $\mathcal T_{C}^{e}(1)$,
where $e=[1,1,\cdots,1]^{\top}$. It is an interesting topic that how
the two schemes (self-organization and driving) affect the cluster
synchronizability for a given graph topology and will be a topic in
the future.

Re-consider the examples in Sec.II.D, we can use Matlab LMI and
Control Toolbox to obtain the numerical values of $CS_{\mathcal
G,\mathcal C}$ for three graphs shown in Fig. \ref{graphs}. Thus, we
can derive their values: $0.472$, $0.178$, and $0.172$,
respectively. Notwithstanding the right-hand of the Lorenz system
does not satisfy the decreasing condition globally, as detailed
analyzed in Ref. \cite{XC}, the trajectory of the coupled Lorenz
systems is essentially bounded, of which the bound is independent of
the coupling strength $c$. So, concentrating on the bounded terminal
region of all trajectories, the decreasing condition can be
satisfied and $\alpha$ can be estimated \footnote{Since these
theoretical estimation is rather loose, we use computer-aided method
to get the estimation}. Here, we get $\alpha\approx 119.3021$,
$120.9882$, and $114.6048$, respectively. Thus, we obtain
estimations of the infimum of $c$: $252.795$ for the graph 1,
$766.892$ for the graph 2, and $667.9655$ for the graph 3. The
details of reasoning and algebras are omitted here. It is clear that
they all locate in the region of cluster synchronization as
numerically illustrated in fig. \ref{var} but less accurate since
the estimations are rather loose.

\subsection{Generalized weighted topologies}
Previous discussions can also be available toward the coupled system
(\ref{Eq.2}) with general weights.
\begin{eqnarray}
\dot{x}^{i}=f_{k}(x^{i})+\sum_{j=1}^{m}a_{ij}w_{ij}\Gamma(
x^{j}-x^{i}),~i\in\mathcal C_{k},~k=\onetoK.\label{Eq.22}
\end{eqnarray}
Here, the graph may be {\bf directed}, i.e., $a_{ij}=1$, if there is
an edge from vertex $j$ to vertex $i$, otherwise, $a_{ij}=0$.
Weights are even not required positive. For the existence of
invariant cluster synchronization manifold, we assume
\begin{eqnarray}
\sum_{j\in\mathcal N_{k'}(i)}w_{ij}=\sum_{j'\in\mathcal
N_{k'}(i')}w_{i'j'}\label{common1}
\end{eqnarray}
holds for all $i,i'\in\mathcal C_{k}$ and $k\ne k'$. Define its
Laplacian $G=[g_{ij}]_{i,j=1}^{m}$ as follows.
\begin{eqnarray*}
g_{ij}=\left\{\begin{array}{ll} w_{ij}&a_{ij}=1\\
0&i\ne j{\rm~and~}a_{ij}=0\\
-\sum\limits_{k=1, k\ne i}^{m}g_{ik}&i=j\end{array}\right..
\end{eqnarray*}
Thus, Eq. (\ref{Eq.22}) becomes
\begin{eqnarray}
\dot{x}^{i}=f_{k}(x^{i})+\sum_{j=1}^{m}g_{ij}\Gamma
x^{j},~i\in\mathcal C_{k},~k=\onetoK.\label{Eq.41}
\end{eqnarray}

Replacing $cl_{ij}$ by $g_{ij}$ and following the routine of the
proof of theorem \ref{st1}, we can obtain following result.
\begin{theorem}\label{thm1}
Suppose that the common inter-cluster coupling condition
(\ref{common1}) is satisfied, each $f_{k}(\cdot)-\alpha\Gamma\cdot$
satisfies the decreasing condition for some $\alpha\in\mathbb R$,
and $\Gamma$ is nonnegative definite. If there exists a positive
definite diagonal matrix $D$ such that
\begin{eqnarray}
\bigg[D(G+\alpha I_{m})\bigg]^{s}_{\mathcal T_{\mathcal
C}^{d}(1)}\le 0\label{gen1}
\end{eqnarray}
holds, then the coupled system (\ref{Eq.41}) can clustering
synchronize with respect to the clustering $\mathcal C$.
\end{theorem}
And, we use the same discussions as in theorem \ref{state2} to
obtain the following general result.
\begin{theorem}\label{thm2}
Suppose that the common inter-cluster coupling condition
(\ref{inter-cluster}) is satisfied. For a bi-directed unweighted
graph $\mathcal G$, there exist positive weights to the graph
$\mathcal G$ such that Ineq. (\ref{gen1}) holds if and only if all
vertices in the same cluster belongs to the same connected component
in the graph $\mathcal G$.
\end{theorem}
In fact, the proofs of theorems \ref{thm1} and \ref{thm2} simply
repeats those of theorems \ref{st1} and \ref{state2}, respectively.

Ref. \cite{Sor} is a paper closely relating to this paper. Here, we
give some comparisons. First, investigated the local cluster
synchronization of inter-connected clusters by extending the master
stability function method. Instead, in this paper, we are concerned
with the global cluster synchronization. Second, the models
discussed are different. The topologies discussed in \cite{Sor}
exclude intra-cluster couplings. In this paper, we consider more
general graph topology. Third, Ref. \cite{Sor} studied the situation
of nonlinear coupling function and we consider the linear case.
Despite that Ref. \cite{Sor} considered different coupling stengths
for clusters and we consider a common one in Sec. II, theorem
\ref{thm1} can apply to discussion of such models proposed in Ref.
\cite{Sor}.

\section{Conclusions}
The idea for studying synchronization in networks of coupled
dynamical systems sheds light on cluster synchronization analysis.
In this paper, we study cluster synchronization in networks of
coupled non-identical dynamical systems. Cluster synchronization
manifold is defined as that the dynamics of the vertices in the same
cluster are identical. The criterion for cluster synchronization is
derived via linear matrix inequality. The differences between
clustered dynamics are guaranteed by the non-identical dynamical
behaviors of different clusters. The algebraic graph theory tells
that the communicability between each pair of vertices in the same
cluster is a doorsill for chaos cluster synchronization. This leads
an description of two schemes to realize cluster synchronization:
self-organization and driving. One can see that the latter scheme
implies that cluster synchronization can be realized in a
non-clustered networks, for example, the graph 2 in the Fig.
\ref{graphs}. Adaptive feedback algorithm is used to enhance cluster
synchronization motions.

\begin{acknowledgments}
This work was jointly supported by the National Natural Sciences
Foundation of China under Grant Nos. 60774074 and , the Mathematical
Tian Yuan Youth Foundation of China No.10826033, and SGST
09DZ2272900.
\end{acknowledgments}

\section*{Appendix}

In this appendix, for each positive $d$, we give the details to find
a $u\in\mathcal T_{\mathcal C}^{d}(1)$ with $u\ne 0$ such that
$u^{\top}DL u=0$ in the case that there exists a cluster $\mathcal
C_{1}$ that does not belong the same connected component. Without
loss of generality, suppose $L$ has form:
\begin{eqnarray*}
L=\left[\begin{array}{ll}L_{1}&0\\
0&L_{2}\end{array}\right].
\end{eqnarray*}
Let $\mathcal V_{1}$ and $\mathcal V_{2}$ corresponds the
sub-matrices $L_{1}$ and $L_{2}$ respectively. And, $\mathcal
C_{1}\bigcap\mathcal V_{i}\ne\emptyset$ for all $i=1,2$. We consider
two situations. First, in the case that $\mathcal C_{1}$ is isolated
from other clusters. In other words, there are no edges between
$\mathcal C_{1}$ and other clusters. Let
\begin{eqnarray*}
u_{i}=\left\{\begin{array}{ll}\alpha&i\in\mathcal
C_{1}\bigcap\mathcal V_{1}\\
\beta&i\in\mathcal C_{1}\bigcap\mathcal V_{2}\\
0&{\rm otherwise}\end{array}\right..
\end{eqnarray*}
Let $a=\sum_{j\in\mathcal C_{1}\bigcap\mathcal V_{1}}d_{j}$ and
$b=\sum_{j\in\mathcal C_{1}\bigcap V_{2}}d_{j}$. Then, if picking
$\alpha$ and $\beta$ satisfying $a\alpha+b\beta=0$ with
$\alpha,\beta\ne 0$, then $u\in\mathcal T_{\mathcal C}^{d}(1)$
holds. In addition, $u^{\top}DLu=0$ due to $Lu=0$.

In the case that $\mathcal C_{1}$ is not isolated, suppose there are
totally $K$ clusters, and $L_1$ and $L_2$ are both connected
(otherwise, we only consider the connection parts of $L_1$ and $L_2$
that contain vertices from $\mathcal{C}_1$), due to the {\em common
inter-cluster coupling condition}, and the absence of isolated
cluster, we have $\mathcal{C}_i\bigcap \mathcal{V}_j\not=\emptyset$
holds for all $i=1,\cdots,K$ and $j=1,2$. Pick a vector
$u=[u_1,\cdots,u_m]^{\top}$ with
\begin{eqnarray*}
u_i=\left\{\begin{array}{ll} \alpha_k& i\in \mathcal{C}_k\bigcap
\mathcal{V}_1\\
\beta_k& i \in \mathcal{C}_k\bigcap \mathcal{V}_2,
\end{array}\right.
\end{eqnarray*}
Denote $d_k^1=\sum_{i\in \mathcal{C}_k\bigcap\mathcal{V}_1}d_i$,
$d_k^2=\sum_{i\in \mathcal{C}_k\bigcap\mathcal{V}_2}d_i$, and
$\bar{u}_1=[\alpha_1,\cdots,\alpha_K]^{\top}$,
$\bar{u}_2=[\beta_1,\cdots,\beta_K]^{\top}$,
$\bar{u}=[\bar{u}_1^{\top},\bar{u}_2^{\top}]^{\top}$,
 $\bar{D}_1={\rm diag}[ d_1^1,\cdots,d_K^1 ]$, $\bar{D}_2={\rm diag}[d_1^2,\cdots,d_K^2]$, and
 $\bar{D}={\rm diag}[ \bar{D}_1,\bar{D}_2 ]$. Define a $K\times K$
 matrix $W^1$ from $L_1$ in such way that for $i\not=j$, $W^1_{ij}=1$ if there's interaction
 between cluster $i$ and $j$, and $W^1_{ij}=0$ otherwise. $W^1_{ii}=-\sum\limits_{j=1,j\not=i}^{K}W^1_{ij}$.
 Define $W^2$ in the same way according to $L_2$, due to the {\it common inter-cluster condition}, it is
 easy to see that $W^1=W^2$. Denote $W={\rm diag}[W^1,W^2]$.

After computation, we have that for any given positive definite
diagonal matrix $D={\rm diag}[d_1,\cdots,d_m]$,
$u^{\top}DLu=\bar{u}^{\top}\bar{D}W\bar{u}$ holds. For $u\in
\mathcal{T}_{\mathcal{C}}^d$,
$\bar{u}_2=-\bar{D}_1\bar{D}_2^{-1}\bar{u}_1$. Denote
$v=\bar{D}_1\bar{u}_1$, we have $
\bar{u}^{\top}\bar{D}W\bar{u}=[v^{\top}
v^{\top}]W\bar{D}^{-1}[v^{\top}
v^{\top}]^{\top}=v^{\top}W^{1}(\bar{D}_{1}^{-1}+\bar{D}_{2}^{-1})v$.
This implies that if we can find $v$ satisfying
$v^{\top}W^1(\bar{D}_1^{-1}+\bar{D}_2^{-1})v=0$, then there exists
$u\in \mathcal T_{\mathcal{C}}^d(1)$ such that $u^{\top}DLu=0$.
Since $W^1(\bar{D}_1^{-1}+\bar{D}_2^{-1})$ has rank at most $K-1$,
we can pick $v$ as the eigenvector corresponding to the zero
eigenvalue of $W^1(\bar{D}_1^{-1}+\bar{D}_2^{-1})$, and this
completes the proof.

In summary, in both situations, we can find certain nonzero vector
$u$ belonging to the transverse space $\mathcal T_{\mathcal
C}^{d}(1)$ and $u^{\top}DLu=0$.



\begin{thebibliography}{99}
\bibitem{interval}
The reason why we choose initial time for average from 50 not 0 is
that even if the coupled system synchronizes clustering, the
variance as calculated by ${\rm var}$ can be very large ($>10^3$).
This implies that it will take a very long time for average to
make the variance is near zero. To save calculationg amount, we
pick the inital time apart from zero.
\bibitem{sth}
Our theoretical result (proposition \ref{st1}) can only intepret
the case the coupled system can clustering synchronize if the
coupling strength $c$ is large enough. As shown by Fig. \ref{var},
the coupled system (\ref{Eq.3}) over the graphs 2 and 3 can
clustering synchronize only if the coupling strength $c$ is
greater than some threshold. But, for the graph 1, despite the
coupled system can synchronize if $c$ is greater than some value
(around 10), there exists an interval (about $[2.2,5]$) of $c$ by
which the system synchronizes. This can not be interpreted via our
theoretical result.
\end{thebibliography}
\end{document}